\documentclass[12pt]{amsart}
\usepackage{amsmath,empheq}
  \usepackage{paralist}
  \usepackage{graphics} 
  \usepackage{epsfig} 
\usepackage{graphicx}  \usepackage{epstopdf}
 \usepackage[colorlinks=true]{hyperref}
 \usepackage{amsfonts}
\usepackage{amssymb}
\usepackage{mathrsfs}
\hypersetup{urlcolor=blue, citecolor=red}
\newtheorem{theorem}{Theorem}[section]

\theoremstyle{definition}

\newcommand{\ep}{\varepsilon}

\newcommand{\RR}{\mathbb R}

\newcommand{\jj}{\mathcal J}
\newcommand{\ee}{\mathcal E}
\newcommand{\ri}{\rightarrow}
\newcommand{\di}{\displaystyle}
\newcommand{\bb}{\begin{equation}}
\newcommand{\bbb}{\end{equation}}
\newcommand{\wrn}{W^{1,p(x)}(\RR^N)}
\newcommand{\wrnr}{W^{1,p(x)}_{{\rm rad}}(\RR^N)}

\title[Infinitely many symmetric solutions] 
      {Infinitely many symmetric solutions for anisotropic problems driven by
       nonhomogeneous operators}
\author[Du\v{s}an D. Repov\v{s}]{}
 \email{dusan.repovs@guest.arnes.si}

\begin{document}
\maketitle
\centerline{\scshape Du\v{s}an D. Repov\v{s}}
\medskip
{\footnotesize
 \centerline{Faculty of Education and Faculty of Mathematics and Physics}
   \centerline{University of Ljubljana}
   \centerline{Institute of Mathematics, Physics and Mechanics}
   \centerline{ SI-1000 Ljubljana, Slovenia}

\bigskip

\centerline{\it Dedicated to Professor Vicen\c{t}iu R\u{a}dulescu for his 60$^{th}$ birthday}

\begin{abstract}
We are concerned with the existence of infinitely many radial symmetric solutions for a nonlinear stationary problem driven by a new class of nonhomogeneous differential operators. Our proof relies on the symmetric version of the mountain pass theorem.
\end{abstract}

\section{Introduction}
Given an even functional on an infinite-dimensional Banach space that fulfills natural assumptions, the symmetric mountain pass lemma of P.~Rabinowitz \cite{rab} establishes the existence of an unbounded sequence of critical values.  This result extends to a symmetric framework the initial version of the mountain pass theorem due to A.~Ambrosetti and P.~Rabinowitz \cite{ambrab}. At the same time, the symmetric mountain pass theorem can be viewed as an extension of the Ljusternik-Schnirelmann theorem in the framework of unbounded functionals defined on Banach spaces. As pointed out by H.~Brezis and
F. Browder \cite{brebro}, the mountain pass theorem ``extends ideas already present at
Poincar\'e and Birkhoff". We refer to Y.~Jabri \cite{jabri} and P.~Pucci and V.~R\u adulescu \cite{pucrad} for more details on the mountain pass theorem and related applications.

We recall  the original statement of the symmetric mountain pass theorem.

\begin{theorem}\label{smp}
Let $X$ be a real infinite-dimensional Banach space and $\jj\in C^1(X,\RR)$ a functional satisfying the Palais-Smale condition and the following hypotheses:

(i) $\jj(0)=0$ and there are constants $\rho,\alpha>0$ such that $\jj_{|\partial B_\rho}\geq\alpha$;

(ii) $\jj$ is even; and

(iii) for all finite-dimensional subspaces $X_0\subset X$, there exists $R=R(X_0)>0$ such that
$$\jj(u)\leq 0\quad\mbox{for all $u\in X_0\setminus B_R(X_0)$}.$$

Then $\jj$ has an unbounded sequence of critical values.
\end{theorem}

This result is an efficient tool for proving multiplicity properties in semilinear or quasilinear elliptic problems with odd nonlinearities and Dirichlet boundary condition. The standard application of Theorem \ref{smp} concerns the following boundary value problem (see Y.~Jabri \cite[pp. 122-124]{jabri})
\begin{equation}\label{yjabri}
\left\{\begin{array}{lll}
&-\Delta u=f(x,u)&\quad\mbox{in}\ \Omega\\
& u=0&\quad\mbox{on}\ \partial\Omega,
\end{array}\right.
\end{equation}
where $\Omega\subset\RR^N$ ($N\geq 3$) is a bounded domain with smooth boundary and $f:\RR\ri\RR$ is a Carath\'eodory function with the following properties:\\
(i) $f$ is odd in $u$, that is, $f(x,-u)=-f(x,u)$;\\
(ii) there exists $p\leq 2^*:=2N/(N-2)$ such that $f$ satisfies the growth condition
$$|f(x,u)|\leq C(1+|u|^{p-1})\quad\mbox{a.e.}\ (x,u)\in\Omega\times\RR;$$
(iii) there are constants $\mu>2$ and $r>0$ such that for almost every $x\in\Omega$ and all $|u|\geq r$
$$0<\mu F(x,u)\leq uf(x,u),\quad\mbox{where}\ F(x,u):=\int_0^uf(x,t)dt.$$

Under these hypotheses, Theorem \ref{smp} yields the existence of an unbounded sequence of solutions of problem \eqref{yjabri}.

The present paper was inspired by recent advances in the study of nonlinear stationary problems driven by nonhomogeneous differential operators. Important pioneering contributions to this field are due to T.C.~Halsey \cite{halsey} and V.V.~Zhikov \cite{zhikov} who studied the behaviour of non-Newtonian electrorheological fluids and anisotropic materials. These models strongly rely on partial differential equations with variable exponent, which have been intensively studied in the last few decades. We refer to the recent monograph by V.~R\u adulescu and D.~Repov\v{s} \cite{radrep} for a comprehensive qualitative analysis of nonlinear PDEs with variable exponent by means of variational and topological methods. These problems (with possible lack of uniform convexity) are essentially described by the differential operator
$$\Delta_{p(x)}u:={\rm div}\,(|\nabla |^{p(x)-2}\nabla u),$$
which changes its growth properties according to the point.
More precisely, the variable exponent $p(x)$ describes
the geometry of a material that is allowed to change its hardening exponent according to the point.
Recently,  I.H. Kim and Y.H. Kim \cite{kim} introduced a new class of nonhomogeneous differential operators, which extend the standard operators with variable exponent. We refer to S.~Baraket, S.~Chebbi, N.~Chorfi, and V.~R\u adulescu \cite{baraket} and N.~Chorfi and V.~R\u adulescu \cite{chorfi} for contributions in this new abstract setting.

In order to introduce the problem studied in this paper and our main result, we need to recall some basic notions and properties. We refer to V.~R\u adulescu and D.~Repov\v{s} \cite{radrep}, resp. I.H.~Kim and Y.H.~Kim \cite{kim} for more details.

\section{Lebesgue and Sobolev spaces with variable exponent}
Let
{\small$$C_+(\RR^N):=\{p:\RR^N\ri\RR; \ p\ \mbox{continuous,}\ 2\leq N<\inf_{x\in\RR^N}p(x)\leq \sup_{x\in\RR^N}p(x)<\infty\}.$$}If $p\in C_+(\RR^N)$, we set $p^-:=\inf_{x\in\RR^N}p(x)$ and $p^+:=\sup_{x\in\RR^N}p(x)$.

For all $p\in C_+(\RR^N)$, let $L^{p(x)}(\RR^N)$ be the Lebesgue space with variable exponent defined by
$$L^{p(x)}(\RR^N):=\left\{u:\RR^N\ri\RR;\ \mbox{$u$ is measurable and}\ \int_{\RR^N}|u(x)|^{p(x)}dx<+\infty\right\}$$
and endowed with the norm
$$|u|_{p(x)}:=\inf\left\{\mu>0;\ \int_{\RR^N}\left|\frac{u(x)}{\mu}\right|^{p(x)}dx\leq 1\right\}.$$

Let $L^{p'(x)}(\RR^N)$ be the dual space
of $L^{p(x)}(\RR^N)$, where $1/p(x)+1/p'(x)=1$. Then for all
$u\in L^{p(x)}(\RR^N)$ and $v\in L^{p'(x)}(\RR^N)$  the following H\"older-type inequality holds:
\begin{equation}\label{Hol}
\left|\int_{\RR^N} uv\;dx\right|\leq\left(\frac{1}{p^-}+
\frac{1}{p'^-}\right)|u|_{p(x)}|v|_{p'(x)}\,.
\end{equation}

Next, we define the corresponding Sobolev function space with variable exponent by
$$W^{1,p(x)}(\RR^N):=\{u\in L^{p(x)}(\RR^N);\ |\nabla u|\in L^{p(x)}(\RR^N)\}.$$
This space is endowed with the norm
$$\|u\|_{p(x)}:=|u|_{p(x)}+|\nabla u|_{p(x)}.$$

The critical Sobolev exponent of $p\in C_+(\RR^N)$ is defined by
$$p^*(x):=\left\{
\begin{array}{lll}
&\di\frac{Np(x)}{N-p(x)}&\quad\mbox{if}\ p(x)<N\\
&\di +\infty&\quad\mbox{if}\ p(x)\geq N.\end{array}\right.
$$

 The function spaces with variable exponent have some striking properties, namely:

(i) If $p\in C_+(\RR^N)$ and $p^+<\infty$, then the formula
$$\int_{\RR^N} |u(x)|^pdx=p\int_0^\infty t^{p-1}\,|\{x\in\RR^N ;\ |u(x)|>t\}|\,dt$$
has  no variable exponent analogue.

(ii) Variable exponent Lebesgue spaces do not have the  mean continuity property. More precisely, if $p$ is continuous and nonconstant in an open ball $B$, then there exists a function $u\in L^{p(x)}(B)$ such that $u(x+h)\not\in L^{p(x)}(B)$ for all $h\in{\mathbb R}^N$ with arbitrary small norm.

(iii) The function spaces with variable exponent
 are {\it never} translation invariant.  The use
of convolution is also limited, for instance  the Young inequality
$$| f*g|_{p(x)}\leq C\, | f|_{p(x)}\, \| g\|_{L^1}$$
holds if and only if
$p$ is constant.

We refer to \cite{radrep} for additional properties.

\section{A new nonhomogeneous differential operator}
Assume that $p\in C_+(\RR^N)$ and consider the mapping $\phi:\RR^N\times [0,\infty)\ri [0,\infty)$ that satisfies the following condiitons:

\smallskip
\noindent (H1) the function $\phi(\cdot,\xi)$ is measurable for all $\xi\geq 0$ and $\phi(x,\cdot)$ is locally absolutely continuous on $[0,\infty)$ for almost all $x\in\Omega$;

\smallskip\noindent (H2) there exist $a\in L^{p'(x)}(\RR^N)$ and $b>0$ such that
$$|\phi (x,|v|)v|\leq a(x)+b|v|^{p(x)-1}$$
for almost all $x\in\Omega$ and for all $v\in\RR^N$;

\smallskip\noindent (H3) there exists $c>0$ such that
$$\phi(x,\xi)\geq c\xi^{p(x)-2},\quad \phi(x,\xi)+\xi\frac{\partial\phi}{\partial\xi}(x,\xi)\geq c\xi^{p(x)-2}$$
for almost all $x\in\Omega$ and for all $\xi>0$.

\smallskip For $\phi$ with the above properties we set
\bb\label{A0def}A_0(x,t):=\int_0^t\phi (x,s)s ds.\bbb
Consider the functional $A:W^{1,p(x)}(\RR^N)\ri\RR$ defined by
$$A(u):=\int_{\RR^N} A_0(x,|\nabla u|)dx.$$

Assume that hypotheses (H1) and (H2) hold. Then by \cite[Lemma 3.2]{kim}, the nonlinear operator $A$ is of class $C^1$ and its G\^ateaux derivative is given by
$$\langle A'(u),v\rangle=\int_{\RR^N}\phi(x,|\nabla u|)\nabla u\nabla vdx,\quad\mbox{for all}\ u,v\in W^{1,p(x)}(\RR^N).$$

Let us now assume that hypotheses (H1)--(H3) are fulfilled. By \cite[Lemma 3.4]{kim}, the operator $A:\wrn\ri W^{1,p'(x)}(\RR^N)$ is strictly monotone and a mapping of type ($S_+$), that is, if $u_n\rightharpoonup u$ in $\wrn$ as $n\ri\infty$ and $\limsup_{n\ri\infty}\langle A'(u_n)-A'(u),u_n-u\rangle\leq 0$, then $u_n\ri u$ in $\wrn$ as $n\ri\infty$.

 The nonhomogeneous differential operator ${\rm div}\, (\phi(x,|\nabla u|)\nabla u)$, where $\phi$ satisfies (H1)--(H3)
was introduced in \cite{kim}.
This operator generalizes the usual operators with variable exponent. For instance, if $\phi (x,\xi)=\xi^{p(x)-2}$ then we obtain the standard $p(x)$-Laplace operator, that is, $\Delta_{p(x)}u:={\rm div}\, (|\nabla u|^{p(x)-2}\nabla u)$.
The new abstract setting includes the case $\phi (x,\xi)=(1+|\xi|^2)^{(p(x)-2)/2}$, which corresponds to the generalized mean curvature operator
$${\rm div}\, \left[(1+|\nabla u|^2)^{(p(x)-2)/2}\nabla u \right].$$
The capillarity equation corresponds to
$$\phi(x,\xi)=\left(1+\frac{\xi^{p(x)}}{\sqrt{1+\xi^{2p(x)}}}\right)\xi^{p(x)-2},\quad x\in\Omega,\ \xi>0,$$ hence the corresponding capillary phenomenon is described by the differential operator
$${\rm div}\, \left[\left( 1+\frac{|\nabla u|^{p(x)}}{\sqrt{1+|\nabla u|^{2p(x)}}} \right)|\nabla u|^{p(x)-2}\nabla u\right].$$

\section{The main result}
Throughout this paper we shall
assume that $p\in C_+(\RR^N)$ and $p$ is a radial function, that is, $p(x)=p(|x|)$ for all $x\in\RR^N$.
Let
$$\wrnr :=\{u\in\wrn ;\ \mbox{$u$ is radial}\}.$$

We are concerned with the study of the following nonlinear problem
\bb\label{problem}
\di -{\rm div}\, (\phi(x,|\nabla u|)\nabla u)+\phi(x,| u|) u= V(x)f(u)\quad\mbox{in}\ \RR^N,
\bbb
where the potential $V:\RR^N\ri [0,+\infty)$ and the nonlinearity $f$ satisfy the following hypotheses:
\bb\label{V}
V\in L^1(\RR^N)\cap L^\infty(\RR^N)\ \mbox{is radial and $\exists\, r_0>0$ such that}\ \inf_{|x|\leq r_0}V(x)>0,
\bbb
\bb\label{f1}
\mbox{$f$ is odd and}\ \lim_{u\ri 0}f(u)/|u|^{p^+-1}=0,
\bbb
and
\bb\label{f2}
\mbox{there exist $\mu>\frac{bp^+}{c}$ and $M>0$ such that $0<\mu F(u)\leq uf(u)$ for all $|u|\geq M$,}
\bbb
where $F(u):=\int_0^uf(t)dt$.

In this paper, due to the symmetry assumptions imposed to $p$ and $V$, we are looking for radial solutions of problem \eqref{problem}.

We say that $u\in\wrnr\setminus\{0\}$ is a solution of problem \eqref{problem} if
$$\int_{\RR^N} \left[\phi(x,|\nabla u|)\nabla u\cdot\nabla v+\phi(x,| u|) u v\right]dx=\int_{\RR^N} V(x)f(u)vdxdx,$$
for all $v\in\wrnr$.

The energy functional associated to problem \eqref{problem} is $\ee:\wrn\ri\RR$ defined by
$$\ee (u):=\int_{\RR^N}\left[A_0(x,|\nabla u|)+A_0(x,|u|) \right]dx-\int_{\RR^N}V(x)F(u)dx.$$

Let $\ee_0$ denote $\ee$ restricted to the function space $\wrnr$. By the isometric Palais principle \cite{palais} (see also \cite[Theorem 1.50]{krv}), any critical point of $\ee_0$ is also a critical point of $\ee$. This shows that finding radially symmetric solutions of problem \eqref{problem} reduces to finding
 the nontrivial critical points of the energy functional $\ee_0$.

\begin{theorem}\label{t1}
Assume that hypotheses (H1)--(H3), \eqref{V}, \eqref{f1}, and \eqref{f2} are fulfilled. Then problem \eqref{problem} has infinitely many solutions.
\end{theorem}

As we shall
see in the proof of this result, problem \eqref{problem} still has at least one (radially symmetric) solution, provided that the oddness symmetry hypothesis on $f$ is removed.

Theorem \ref{t1} extends the pioneering multiplicity result of A.~Ambrosetti and P.~Rabinowitz \cite[Theorem 3.13]{ambrab} in the following directions:

(i) the standard (linear) second order uniformly elliptic operator $\di\Sigma_{i,j=1}^N(a_{ij}(x)\break u_{x_i})_{x_j}$ is replaced by the nonhomogeneous differential operator ${\rm div}\, (\phi(x,|\nabla u|)\nabla u)$;

(ii) our study is performed in the entire Euclidean space, instead of a bounded domain with smooth boundary. However, in our abstract setting, the lack of compactness of $\RR^N$ is compensated by the compactness of the embedding of the space $\wrnr$ into $L^\infty(\RR^N)$, provided that $N<p^-\leq p^+<+\infty$.
	
\section{Proof of Theorem \ref{t1}} We first check that $\ee_0$ satisfies
the geometric hypotheses of the mountain pass theorem.

\smallskip\noindent{\bf Step 1.} There exist positive  constants $\rho$ and $\alpha$ such that $\ee_0 (u) \geq\alpha$ for all $u\in \wrnr$ with $\|u\|=\rho$.

For $\rho\in (0,1)$ (to be prescribed later), we fix $u\in\wrnr$ satisfying $\|u\|_{p(x)}=\rho$. By hypothesis (H3),
 we have
$$\int_{\RR^N}\left[A_0(x,|\nabla u|)+A_0(x,|u|) \right]dx \geq \frac{c}{p^+}\int_{\RR^N}(|\nabla u|^{p(x)}+|u|^{p(x)})dx.$$
Next, we use relation (1.8) in \cite[p. 11]{radrep}. Thus, since $\|u\|_{p(x)}=\rho<1$, we have
\bb\label{1eq}
\int_{\RR^N}\left[A_0(x,|\nabla u|)+A_0(x,|u|) \right]dx \geq \frac{c}{p^+}\,\|u\|^{p^+}_{p(x)}.
\bbb

On the other hand, assumption \eqref{f1} implies that $F(u)/|u|^{p^+}\ri 0$ as $u\ri 0$. Fix $\ep>0$. It follows that if $\rho>0$ is small enough then
$$\int_{\RR^N}V(x)F(u)dx\leq\ep\,\|V\|_{L^1}\, \|u\|^{p^+}_{L^\infty}.$$
Since $\wrnr$ is continuously embedded into $L^\infty(\RR^N)$, we deduce that there exists $C>0$ such that for all $u\in\wrnr$ with $\|u\|_{p(x)}=\rho$ we have
\bb\label{2eq}
\int_{\RR^N}V(x)F(u)\leq C\ep\,\|u\|^{p^+}_{p(x)}.\bbb

Combining relations \eqref{1eq} and \eqref{2eq},
 we obtain
$$\begin{array}{ll}
\ee_0(u)&\di\geq \frac{c}{p^+}\,\|u\|^{p^+}_{p(x)}-C\ep\,\|u\|^{p^+}_{p(x)}\\
&=\frac{c}{p^+}\,\rho^{p^+}-C\ep\, \rho^{p^+}.
\end{array}$$
Choosing $\ep=c/(2Cp^+)>0$, we have
$$\ee_0(u)\geq \frac{c}{2p^+}\,\rho^{p^+}=:\alpha>0.$$

\smallskip\noindent{\bf Step 2.}
For all finite-dimensional subspaces $X_0\subset \wrnr$
 there exists $R=R(X_0)>0$ such that
$$\ee_0(u)\leq 0\quad\mbox{for all $u\in X_0\setminus B_R(X_0)$}.$$

We first claim that for all $w\in\wrnr$ with $\|w\|_{p(x)}=1$ there exists $\lambda (w)>0$ such that
\bb\label{claim}\ee_0(\lambda w)<0\ \mbox{for all $\lambda\in\RR$ with $|\lambda|\geq \lambda (w)$}.\bbb

We observe that hypothesis \eqref{f2} implies that there are positive constants $C_1$ and $C_2$ such that
$$f(u)\geq C_1u^{\mu-1}-C_2\quad\mbox{for all $u\geq 0$.}$$
Therefore
$$F(u)\geq \frac{C_1}{\mu}\, u^{\mu}-C_2u\quad\mbox{for all $u\geq 0$.}$$
It follows that there exists $C_0>0$ such that
\bb\label{Fmare}
F(u)\geq C_0 u^{\mu}\quad\mbox{for all $|u|\geq M$.}\bbb

Fix $w\in\wrnr$ with $\|w\|_{p(x)}=1$ and $\lambda\in\RR$ (with $|\lambda|>1$).

In order to find an upper estimate for $\ee_0(\lambda w)$, we first observe that using hypothesis (H2), we have
\bb\label{esti1}\begin{array}{ll}
&\di \int_{\RR^N}[A_0(x,|\lambda \nabla w)+A_0(x,|\lambda w|)]dx\di \leq \\
&\di |\lambda|\int_{\RR^N}a(x)(|\nabla w|+|w|)dx+b\int_{\RR^N}
\left(\int_0^{|\lambda \nabla w|}s^{p(x)-1}ds+\int_0^{|\lambda  w|}s^{p(x)-1}ds \right)dx\leq\\
&\di  |\lambda|\int_{\RR^N}a(x)(|\nabla w|+|w|)dx+\frac{b|\lambda|^{p^+}}{p^-}\,\int_{\RR^N}(|\nabla w|^{p(x)}+|w|^{p(x)})dx=\\
&\di \frac{b|\lambda|^{p^+}}{p^-}+C|\lambda|,
\end{array}
\bbb
where $C>0$ is a constant depending only on $|a|_{p'(x)}$ and the best constant of the continuous embedding $\wrnr\hookrightarrow L^{p(x)}(\RR^N)$.

By \eqref{Fmare} we have
\bb\label{esti2}
\int_{|\lambda w|\geq M}V(x)F(\lambda w)dx\geq C_0\int_{|\lambda w|\geq M}|\lambda|^\mu |w|^\mu dx=C_0|\lambda|^\mu \int_{|\lambda w|\geq M}V(x)|w|^\mu dx.\bbb

Since $F$ is bounded on the interval $[-M,M]$, there exists $C>0$ such that $F(t)\geq -C$ for all $t\in [-M,M]$. It follows that
\bb\label{esti3}
\int_{|\lambda w|< M}V(x)F(\lambda w)dx\geq -C\int_{|\lambda w|<M}V(x)dx\geq -C\int_{\RR^N}V(x)dx=-C\|V\|_{L^1}\,.\bbb

Estimates \eqref{esti1}, \eqref{esti2}, and \eqref{esti3} imply
$$\ee_0(\lambda w)\leq \frac{b}{p^-}\, |\lambda|^{p^+}+C|\lambda|-C_0|\lambda|^\mu \int_{|\lambda w|\geq M}V(x)|w|^\mu dx+C\|V\|_{L^1},$$
hence $\ee_0(\lambda w)\ri-\infty$ as $|\lambda|\ri\infty$. This proves our claim \eqref{claim}.

Next, since the space $\wrnr$ is compactly embedded into $L^\infty(\RR^N)$, we deduce that there exists $C>0$ such that
$$|w(x)|\leq C\ \mbox{for all $w\in\wrnr$ with $\|w\|_{p(x)}=1$}.$$

This fact implies that our initial claim \eqref{claim} can be improved as follows: for all $w\in\wrnr$ with $\|w\|_{p(x)}=1$ there exist $\lambda (w)>0$ and $\eta(w)>0$ such that
\bb\label{claim2}\ee_0(\lambda z)<0\quad \mbox{$\forall\,|\lambda|\geq \lambda (w)$, $\forall\,z\in\wrnr$, $\|z\|_{p(x)}=1$, $\|z-w\|\leq\eta(w)$}.\bbb

Returning to Step 2, let $X_0\subset \wrnr$ be a finite-dimensional subspace.
Thus the set $X_0\cap\{w\in\wrnr ;\ \|w\|_{p(x)}=1\}$ is compact. Next, using \eqref{claim2}, we deduce that there exists $\lambda_0>0$ depending only on $X_0$ such that
$$\ee_0(\lambda w)\leq 0\quad\mbox{for all $|\lambda|\geq\lambda_0$ and for all $w\in X_0$, $\|w\|_{p(x)}=1$}.$$
Choosing $R(X_0)=\lambda_0$, we obtain the statement contained in our Step 2.

\medskip\noindent{\bf Step 3.} Any Palais-Smale sequence of $\ee_0$ is bounded.

 We  recall that $(u_n)\subset\wrnr$ is a Palais-Smale sequence of $\ee_0$ if
\bb\label{psm} \ee_0(u_n)=O(1)\quad\mbox{and}\quad \|\ee_0'(u_n)\|=o(1)\quad \mbox{as $n\ri\infty$}.\bbb
We also recall that $\ee_0$ satisfies the Palais-Smale condition if any Palais-Smale sequence $(u_n)\subset\wrnr$ of $\ee_0$ is relatively compact.

Arguing by contradiction and normalizing, we assume that (up to a subsequence) $u_n=\lambda_nv_n$, where $\lambda_n=\|u_n\|_{p(x)}\ri\infty$ and $\|v_n\|_{p(x)}=1$. By \eqref{psm} we deduce that
\bb\label{ps0}
\begin{array}{ll}
o(\|u_n\|)&\di =\langle\ee_0'(u_n),u_n\rangle=\int_{\RR^N}[\phi (x,|\nabla u_n|)|\nabla u_n|^2+\phi (x,| u_n|) u_n^2]dx\\
&\di -\int_{[\lambda_n|v_n|\geq M]}V(x)f(\lambda_nv_n)\lambda_nv_ndx-\int_{[\lambda_n|v_n|< M]}V(x)f(\lambda_nv_n)\lambda_nv_ndx.
\end{array}\bbb
Using hypothesis (H2), relation \eqref{ps0} yields
$$
\begin{array}{ll}\lambda_no(1)&\di\leq \lambda_n\int_{\RR^N}a(x)(|\nabla v_n|+|v_n|)dx+b\int_{\RR^N}\lambda_n^{p(x)}(|\nabla v_n|^{p(x)}+|v_n|^{p(x)})dx\\
&\di -\int_{[\lambda_n|v_n|\geq M]}V(x)f(\lambda_nv_n)\lambda_nv_ndx-\int_{[\lambda_n|v_n|< M]}V(x)f(\lambda_nv_n)\lambda_nv_ndx.
\end{array}$$
It follows that
\bb\label{ps1}
\begin{array}{ll} &\di \int_{[\lambda_n|v_n|\geq M]}V(x)f(\lambda_nv_n)\lambda_nv_ndx \di\leq \lambda_n\int_{\RR^N}a(x)(|\nabla v_n|+|v_n|)dx+\\ &\di b\int_{\RR^N}\lambda_n^{p(x)}(|\nabla v_n|^{p(x)}+|v_n|^{p(x)})dx-\int_{[\lambda_n|v_n|< M]}V(x)f(\lambda_nv_n)\lambda_nv_ndx+\lambda_no(1).
\end{array}\bbb
By hypothesis \eqref{f2}, we deduce that
\bb\label{ps2}
\int_{[\lambda_n|v_n|\geq M]}V(x)F(\lambda_nv_n)dx \leq\frac 1\mu \int_{[\lambda_n|v_n|\geq M]}V(x)f(\lambda_nv_n)\lambda_nv_ndx .\bbb

Using hypothesis (H3), we have
\bb\label{ps3}
\begin{array}{ll}
&\di \int_{\RR^N}\left[A_0(x,|\nabla u_n|)+A_0(x,|u_n|) \right]dx=\\ &\di
\int_{\RR^N}\left( \int_0^{|\nabla u_n|}s\phi (x,s)ds+\int_0^{| u_n|}s\phi (x,s)ds\right)dx\geq\\
&\di c\int_{\RR^N}\frac{1}{p(x)}(|\nabla u_n|^{p(x)}+|u_n|^{p(x)})dx=c\int_{\RR^N}\frac{\lambda_n^{p(x)}}{p(x)}\,(|\nabla v_n|^{p(x)}+|v_n|^{p(x)})dx\geq\\
&\di \frac{c}{p^+}\int_{\RR^N}\lambda_n^{p(x)}\,(|\nabla v_n|^{p(x)}+|v_n|^{p(x)})dx.
\end{array}\bbb

Combining relations \eqref{ps2} and \eqref{ps3} in relationship with \eqref{psm}, we deduce
$$\begin{array}{ll}
O(1)&\di =\ee_0(u_n)\geq \frac{c}{p^+}\int_{\RR^N}\lambda_n^{p(x)}\,(|\nabla v_n|^{p(x)}+|v_n|^{p(x)})dx\\
&\di -\frac 1\mu \int_{[\lambda_n|v_n|\geq M]}V(x)f(\lambda_nv_n)\lambda_nv_ndx -\int_{[\lambda_n|v_n|< M]}V(x)F(\lambda_nv_n)dx.\end{array}$$
Using now the estimate established in \eqref{ps1}, we obtain
$$
\begin{array}{ll}
&O(1)\di \geq \frac{c}{p^+}\int_{\RR^N}\lambda_n^{p(x)}\,(|\nabla v_n|^{p(x)}+|v_n|^{p(x)})dx-
\frac{\lambda_n}{\mu}\int_{\RR^N}a(x)(|\nabla v_n|+|v_n|)dx\\
&-\di \frac{b}{\mu}\int_{\RR^N}\lambda_n^{p(x)}(|\nabla v_n|^{p(x)}+|v_n|^{p(x)})dx\\ &\di +\frac{1}{\mu}\int_{[\lambda_n|v_n|< M]}V(x)f(\lambda_nv_n)\lambda_nv_ndx-\int_{[\lambda_n|v_n|< M]}V(x)F(\lambda_nv_n)dx+\lambda_no(1)\\
&\di =\left(\frac{c}{p^+}-\frac b\mu\right)\int_{\RR^N}\lambda_n^{p(x)}(|\nabla v_n|^{p(x)}+|v_n|^{p(x)})dx-\frac{\lambda_n}{\mu}\int_{\RR^N}a(x)(|\nabla v_n|+|v_n|)dx\\
&\di + \frac{1}{\mu}\int_{[\lambda_n|v_n|< M]}V(x)f(\lambda_nv_n)\lambda_nv_ndx-\int_{[\lambda_n|v_n|< M]}V(x)F(\lambda_nv_n)dx+\lambda_no(1).
\end{array}$$
Therefore
\bb\label{ps4}\begin{array}{ll}
O(1)&\di =\ee_0(u_n)\geq \left(\frac{c}{p^+}-\frac b\mu\right)\int_{\RR^N}\lambda_n^{p(x)}(|\nabla v_n|^{p(x)}+|v_n|^{p(x)})dx\\ &\di-\frac{\lambda_n}{\mu}\int_{\RR^N}a(x)(|\nabla v_n|+|v_n|)dx+\lambda_no(1).\end{array}\bbb

Combining hypothesis \eqref{f2} with relation \eqref{ps4} and the fact that $\lim_{n\ri\infty}\lambda_n=+\infty$, we deduce that $\ee_0(u_n)\ri+\infty$ as $n\ri\infty$, which contradicts \eqref{psm}. We conclude that the Palais-Smale sequence $(u_n)$ is bounded in $\wrnr$. We shall
 prove that
  in fact, the sequence converges strongly (up to a subsequence) in $\wrnr$. For this purpose, we apply some ideas developed by R.~Filippucci, P.~Pucci and V.~R\u adulescu \cite[pp. 712-713]{fpr}.

\medskip\noindent
{\bf Step 4.} Any Palais-Smale sequence of $\ee_0$ is relatively compact in $\wrnr$.

Let $(u_n)\subset\wrnr$ be an arbitrary Palais-Smale sequence of $\ee_0$. By Step 3, $(u_n)$ is bounded. Using the compact embedding $\wrnr\hookrightarrow L^\infty(\RR^N)$,
 we can assume that, up to a subsequence,
\bb\label{weak}u_n\rightharpoonup u\quad\mbox{in $\wrnr$}\bbb
\bb\label{strong}u_n\ri u\quad\mbox{in $L^\infty(\RR^N)$.}\bbb
By hypotheses (H2) and (H3) in conjunction with relations \eqref{weak} and \eqref{strong}, we deduce that
\bb\label{23}\int_{\RR^N}\phi(x,|u_n|)u_n(u_n-u)dx-\int_{\RR^N}\phi(x,|u|)u(u_n-u)dx\ri 0\quad\mbox{as}\ n\ri\infty\bbb
and
\bb\label{24}\int_{\RR^N}\phi(x,|\nabla u_n|)\nabla u_n\cdot\nabla(u_n-u)dx-\int_{\RR^N}\phi(x,|\nabla u|)\nabla u\cdot\nabla(u_n-u)dx\ri 0\quad\mbox{as}\ n\ri\infty.\bbb
Recall that Proposition 3.3 in \cite{kim} established that under hypotheses (H1) and (H3) the following Simon-type inequality holds:
$$\phi (x,|u|)u(u-v)-\phi (x,|v|)v(u-v)\geq 4^{1-p^+}c\,|u-v|^{p(x)}\quad\mbox{for all}\ x\in\Omega.$$
Using this inequality, relations \eqref{23} and \eqref{24} imply that
$$u_n\ri u\quad\mbox{in $\wrnr$}.$$
We conclude that any Palais-Smale sequence of $\ee_0$ is relatively compact in\break $\wrnr$, hence $\ee_0$ satisfies the Palais-Smale condition.

Next, we observe that $\ee_0$ is even and $\ee_0(0)=0$. The above steps show that $\ee_0$ fulfills the hypotheses of Theorem \ref{smp}. We deduce that problem \eqref{problem} has infinitely many solutions in $\wrnr$. \qed

\section{Epilogue} The methods developed in this paper can be applied to other recent classes of nonhomogeneous differential operators. For instance, the operator ${\rm div}\, (\phi(x,|\nabla u|)\nabla u)$ that describes problem \eqref{problem} can be replaced with
${\rm div}\, [\phi'(|\nabla u|^2)\nabla u]$, where
for some $1<p<q<N$, the function $\phi\in C^1(\RR_+,\RR_+)$ satisfies the following conditions:

\smallskip
\noindent ($\phi_1$) $\phi (0)=0;$

\smallskip
\noindent ($\phi_2$) there exists $c_1>0$ such that $\phi(t)\geq c_1t^{p/2}$ if $t\geq 1$ and   $\phi(t)\geq c_1t^{q/2}$ if $0\leq t\leq 1$;

\smallskip
\noindent ($\phi_3$) there exists $c_2>0$ such that $\phi(t)\leq c_2t^{p/2}$ if $t\geq 1$ and   $\phi(t)\leq c_2t^{q/2}$ if $0\leq t\leq 1$;

\smallskip
\noindent ($\phi_4$) there exists $0<\mu<1/s$ such that $2t\phi'(t)\leq s\mu\phi (t)$ for all $t\geq 0$;

\smallskip
\noindent ($\phi_5$) the mapping $t\mapsto \phi (t^2)$ is strictly convex.

\smallskip This operator was introduced by A.~Azzollini, P.~d'Avenia, and A.~Pomponio \cite{azo2} and it is described by a potential with different growth near zero and at infinity (the double-power growth hypotheses). We refer to A.~Azzollini \cite{azo1} and  N.~Chorfi and V. R\u adulescu \cite{crejqtde} for recent contributions in connection with the abstract setting generated by this operator.

In a general framework, the presence of two variable exponents $p_1(x)$ and $p_2(x)$ dictates
 the geometry of a composite that changes its hardening  exponent according to the point. Problems with nonstandard growth conditions of $(p,q)$-type have been initially studied by P.~Marcellini \cite{marce} who was interested in the properties of the integral energy functional
$\int_{\RR^N} F(x,\nabla u)dx,$
where $F:\RR^N\times\RR^N\ri\RR$ satisfies unbalanced polynomial growth conditions, namely
$$|\xi|^p \lesssim F(x,\xi)\lesssim |\xi|^q\quad\mbox{with $1<p<q$.}$$

We believe that the main result of this paper can be extended to ``unbalanced" anisotropic differential operators of the type
$$-{\rm div}\, (\phi(x,|\nabla u|)\nabla u)-{\rm div}\, (a(x)\psi(x,|\nabla u|)\nabla u)$$
and
$$-{\rm div}\, (\phi(x,|\nabla u|)\nabla u)-{\rm div}\, (a(x)\psi(x,|\nabla u|)\log (e+|x|)\nabla u).$$

 This abstract setting is in close relationship with the recent contributions of G.~Mingione {\it et al.} \cite{mingi1, mingi2}, who studied non-autonomous problems with associated energies of the type
$$ \int_{\RR^N} [|\nabla u|^{p_1(x)}+a(x)|\nabla u|^{p_2(x)}]dx$$
and
$$ \int_{\RR^N} [|\nabla u|^{p_1(x)}+a(x)|\nabla u|^{p_2(x)}\log(e+|x|)]dx,$$
where $p_1(x)\leq p_2(x)$, $p_1\not=p_2$, and $a(x)\geq 0$.

\section*{Acknowledgements}{This work was supported by the Slovenian Research Agency grants P1-0292, J1-8131, J1-7025, and N1-0064.}

\end{document}